# Tail probabilities for infinite series of regularly varying random vectors

HENRIK HULT[1] and GENNADY SAMORODNITSKY[2]

[1]*Division of Applied Mathematics, Brown University, 182 George Street, Providence, RI 02912, USA.* E-mail: *henrik_hult@brown.edu*
[2]*School of ORIE, Cornell University, 220 Rhodes Hall, Ithaca, NY 14853, USA.*
E-mail: *gennady@orie.cornell.edu*

A random vector $X$ with representation $X = \sum_{j \geq 0} A_j Z_j$ is considered. Here, $(Z_j)$ is a sequence of independent and identically distributed random vectors and $(A_j)$ is a sequence of random matrices, 'predictable' with respect to the sequence $(Z_j)$. The distribution of $Z_1$ is assumed to be multivariate regular varying. Moment conditions on the matrices $(A_j)$ are determined under which the distribution of $X$ is regularly varying and, in fact, 'inherits' its regular variation from that of the $(Z_j)$'s. We compute the associated limiting measure. Examples include linear processes, random coefficient linear processes such as stochastic recurrence equations, random sums and stochastic integrals.

## 1. Introduction

A general and useful class of stochastic models is the class of random coefficient linear models. This is the class of $d$-dimensional random vectors with the stochastic representation

$$X = \sum_{j=0}^{\infty} A_j Z_j. \tag{1.1}$$

Here, $(Z_j)_{j \geq 0}$ is a sequence of independent and identically distributed (i.i.d.) random vectors in $\mathbb{R}^p$, which will be assumed to have a regularly varying law; a precise definition will be given below. A generic element of this sequence is denoted by $Z$. The sequence $(A_j)_{j \geq 0}$ consists of random matrices. Several examples of stochastic models with representation (1.1) are given in Section 2. These include linear processes, solutions to stochastic recurrence equations and stochastic integrals. For the latter two examples, it would be very restrictive to require independence between the coefficients $(A_j)$ and the noise variables $(Z_j)$. Instead, it will be assumed that the coefficients satisfy a certain 'predictability' assumption with respect to the noise variables. The precise formulation of this assumption is given in Section 2.

Naturally, the infinite series in (1.1) is assumed to converge almost surely. Sufficient conditions for such convergence will be provided in the sequel (Theorem 3.1). In this paper, we are interested in the asymptotic tail behavior of $X$. In particular, we study situations under which the random matrices $(A_j)$ are 'small' relative to the noise vectors $(Z_j)$, and the tail probabilities of the random vector $X$ in (1.1) are 'inherited' from the tail probabilities of the noise vectors.

The interest in models with representation (1.1) where the $(Z_j)$ have regularly varying law is motivated by applications. Probability distributions with regularly varying tails have become important building blocks in a wide variety of stochastic models. Evidence for power-tail distributions is by now well documented in a large number of applications including computer networks, telecommunications, finance, insurance, hydrology, atmospheric sciences, geology, ecology, etc. Many of these applications are genuinely multidimensional, involving, for example, multiple servers in a computer network, portfolios of stocks or insurance, or measurements taken at multiple geographical locations. Hence, our interest in the multidimensional case. This requires using the notion of multivariate regular variation, which we give below. Interestingly, once the right notion of regular variation is used, the added generality of the multidimensional result does not require significant additional effort to be established. Indeed, most of the effort in proving the main theorem (Theorem 3.1) is spent on proving the tight moment conditions on the coefficients $(A_j)$.

The notion of multivariate regular variation we will use in this paper is as follows. We say that a $d$-dimensional random vector $Z$ has a regularly varying distribution if there exists a non-null Radon measure $\mu$ on $\overline{\mathbb{R}}^d \setminus \{0\}$ (where $\overline{\mathbb{R}}^d = [-\infty, \infty]^d$), with $\mu(\overline{\mathbb{R}}^d \setminus \mathbb{R}^d) = 0$, such that

$$\frac{P(u^{-1}Z \in \cdot)}{P(|Z| > u)} \xrightarrow{v} \mu(\cdot), \quad (1.2)$$

as $u \to \infty$, on $\overline{\mathbb{R}}^d \setminus \{0\}$. Here, $\xrightarrow{v}$ denotes vague convergence (see Kallenberg (1983), Resnick (1987, 2006)). The limiting measure $\mu$ necessarily obeys a homogeneity property: there is an $\alpha > 0$ such that $\mu(uB) = u^{-\alpha}\mu(B)$ for all Borel sets $B \subset \overline{\mathbb{R}}^d \setminus \{0\}$. This follows from standard regular variation arguments (see, e.g., Resnick (2006), Hult and Lindskog (2006), Theorem 3.1). We will write $Z \in \mathrm{RV}(\mu, \alpha)$ for a random vector satisfying (1.2). For more on multivariate regular variation, we refer to Basrak (2000) and Resnick (1987, 2006).

The tail behavior of $X$ has been previously studied, in particular, in the one-dimensional case $p = d = 1$. In that case, the most general results have been obtained by Resnick and Willekens (1991) and Wang and Tang (2006) under the assumption that the sequence $(A_j)$ is non-negative and independent of the sequence $(Z_j)$. Even in that particular case, our conditions in the case $\alpha \geq 1$ are strictly weaker (the conditions are identical for $0 < \alpha < 1$).

Before we proceed, we would like to point out why the 'predictability' assumption is important. Note that each term in the sum (1.1) is of the form $A_j Z_j$ and $Z_j$ has regularly varying law. If $A_j$ and $Z_j$ are independent and $A_j$ satisfies a moment condition, then a multidimensional version of Breiman's theorem (e.g., Basrak *et al.* (2002), Proposition A.1) can be applied to obtain the tail behavior of the product $A_j Z_j$. As we will see, it is not necessary to assume independence of the entire sequences $(A_j)$ and $(Z_j)$. The 'predictability' assumption which guarantees that $A_j$ and $Z_j$ are independent for each $j$ is sufficient. For a trivial example of how our results can fail without such an assumption, consider the one-dimensional case with $A_0 = I(|Z_0| \leq 1)$ and $A_j = 0$ for $j \geq 1$.

In the next section, we specify the mentioned 'predictability' assumption and provide a wide range of examples that are covered by the general representation (1.1). In Section 3, we state the main theorem of this paper, giving sufficient conditions for convergence of the series $X$ in (1.1) and for $X$ acquiring the regular variation properties from the noise vectors $(Z_j)$. We then explore the implications in the various examples specified in Section 2. The assumptions in the main theorem turn out to be very tight and improve the existing results in most special cases we are considering. The proof of the main theorem is given in Section 4.

## 2. Examples

The representation (1.1) covers many interesting examples of stochastic models that are widely used in applications. Some examples are presented below. We informally divide them into three main categories: 'linear processes', 'stochastic recurrence equations' and 'stochastic integrals'. For the linear processes, it is quite natural to assume independence between the sequence of coefficients $(A_j)$ and the noise sequence $(Z_j)$. However, for the stationary solution to a stochastic recurrence equation one does not, in general, have independence between the coefficients and the noise terms; see Section 2.2 below. There, a 'predictability' assumption is required. For stochastic integrals, it is natural to think of $(A_j)$ as the integrand and introduce a 'predictability' assumption with respect to the sequence $(Z_j)$. To cover all these cases, we introduce the following assumption on the sequences $(A_j)$ and $(Z_j)$.

Assume that there is a filtration $(\mathcal{F}_j, j \geq 0)$ such that

$$A_j \in \mathcal{F}_j, \qquad Z_j \in \mathcal{F}_{j+1}, \qquad \text{for } j \geq 0, \tag{2.1}$$

$$\mathcal{F}_j \text{ is independent of } \sigma(Z_j, Z_{j+1}, \ldots) \qquad \text{for } j \geq 0. \tag{2.2}$$

In a sense, the sequence $(A_j)$ is predictable with respect to the sequence $(Z_j)$. The case where the sequence $(A_j)$ is independent of the sequence $(Z_j)$ is covered by the 'predictable' framework by setting $\mathcal{F}_j = \sigma((A_k)_{k \geq 0}, Z_0, \ldots, Z_{j-1})$.

### 2.1. Linear processes

In this section, we provide detailed examples in the 'linear processes' category.

***Example 2.1 (Linear process).*** Let $(A_j)_{j\geq 0}$ be a sequence of deterministic real-valued $d \times p$ matrices. Then, assuming convergence, $X_k = \sum_{j\geq 0} A_j W_{k-j}$ is a linear process. It is, clearly, stationary. The ($d$-dimensional) distribution of $X_0$ has the representation (1.1) with $Z_j = W_{-j}$.

***Example 2.2 (Random coefficient linear process).*** This is a generalization of the linear process. Let $(\mathbf{A}_k, W_k)_{k\in\mathbb{Z}}$ be a stationary sequence, where each $\mathbf{A}_k = (A_{k,j})_{j\geq 0}$ is itself a sequence of random $d \times p$ matrices. Assuming, once again, convergence, the process $X_k = \sum_{j\geq 0} A_{k,j} W_{k-j}$ is a random coefficient linear process. The stationarity of the sequence $(\mathbf{A}_k, W_k)_{k\in\mathbb{Z}}$ implies that it is a stationary process and the distribution of $X_0$ has the representation (1.1) with $Z_j = W_{-j}$. To apply our results, we will need to assume that the sequences $(A_{0,i})_{i\geq 0}$ and $(W_k)_{k\leq 0}$ satisfy the predictability assumptions (2.1) and (2.2). Of course, assuming that the sequences $(A_{0,i})_{i\geq 0}$ and $(W_k)_{k\leq 0}$ are independent is sufficient for this purpose and this is often the case of interest in the context of linear processes. However, certain important cases of random coefficient linear process are naturally stated in the predictable framework; see Section 2.2.

***Example 2.3 (Partial sum process).*** Consider the random coefficient linear process $(X_k)$ in the previous example and let $S_n = X_1 + \cdots + X_n$. Then, with $B_{n,i} = \sum_{k=1\vee i}^{n} A_{k,k-i}$,

$$S_n = \sum_{k=1}^{n} \sum_{j\geq 0} A_{k,j} W_{k-j} = \sum_{i=-\infty}^{n} \sum_{k=1\vee i}^{n} A_{k,k-i} W_i = \sum_{j\geq 0} B_{n,n-j} W_{n-j}.$$

Of course, $(S_n)_{n\geq 1}$ is not, in general, a stationary process, but at each time $n$, its marginal distribution has the representation (1.1) with $A_j = B_{n,n-j}$ and $Z_j = W_{n-j}$. Once again, predictability assumptions must hold. Simply assuming that the sequence $(\mathbf{A}_k)_{k\geq 0}$ is independent of the sequence $(W_k)_{k\geq 0}$ is sufficient, but not necessary, for this purpose.

## 2.2. Stochastic recurrence equations

A very important particular case of the random coefficient linear model is the stationary solution of a stochastic recurrence equation (SRE).

Assume that $p = d$ and let $(M_k, W_k)_{k\in\mathbb{Z}}$ be a sequence of independent and identically distributed pairs of $d \times d$ matrices and $d$-dimensional random vectors. Under the assumptions $\mathrm{E}\log^+ \|M_k\| < \infty$ and $\mathrm{E}\log^+ |W_k| < \infty$, there exists a strictly stationary causal solution to the SRE

$$X_k = M_k X_{k-1} + W_k, \qquad k \in \mathbb{Z},$$

if and only if the top Lyapunov exponent

$$\gamma = \inf_{n\in\mathbb{N}} \mathrm{E}[(n+1)^{-1} \log \|M_0\| \cdots \|M_{-n}\|] < 0;$$

see Bougerol and Picard (1992), Theorem 1.1, attributed to Brandt (1986). A sufficient condition for $\gamma < 0$ is that $\mathrm{E}\log^+ \|M_0\| < 0$. The stationary solution can be represented as

$$X_0 = W_0 + \sum_{j \geq 1} M_0 \cdots M_{-j+1} W_{-j},$$

which is a random coefficient linear model of the form (1.1) with $A_0 = 1$, $A_j = M_0 \cdots M_{-j+1}, j \geq 1$, and $Z_j = W_{-j}$ (see, e.g., Kesten (1973), Konstantinides and Mikosch (2005)). Note that with $\mathcal{F}_j = \sigma(M_0, W_0, \ldots, M_{-j+1}, W_{-j+1})$, it follows that $A_j \in \mathcal{F}_j$, $Z_j \in \mathcal{F}_{j+1}$ and $\mathcal{F}_j$ is independent of $\sigma(Z_j, Z_{j+1}, \ldots)$. Hence, the predictability assumptions (2.1)–(2.2) are satisfied.

### 2.3. Stochastic integrals

In this section, we provide examples within the 'stochastic integrals' category.

*Example 2.4 (Random sum).* Let $p = d$ and $A_j = \mathrm{Id}\, I\{1 \leq j \leq N\}$, where $N$ is a positive integer-valued random variable which is a stopping time with respect to a filtration $(\mathcal{G}_j, j \geq 0)$ to which the sequence $(Z_j)$ is adapted, and such that $\mathcal{G}_j$ is independent of $\sigma(Z_{j+1}, Z_{j+2}, \ldots)$ for $j \geq 0$. Then $X = \sum_{j=1}^N Z_j$ is a random sum.

*Example 2.5 (Stochastic integral).* The following is a modification of the random sum in the previous example. Let $(C_t)_{t \geq 0}$ be a $p$-dimensional renewal reward process with renewal times $\tau_1, \tau_2, \ldots$. We take a version of $(C_t)$ with right-continuous paths with left limits. For some fixed time $T > 0$, let $N_T$ be the number of renewals until time $T$. Suppose the jump distribution of $(C_t)$ (e.g., the reward distribution) is regularly varying, that is, $\Delta C_{\tau_j} = (C_{\tau_j} - C_{\tau_j-}) \stackrel{\mathrm{d}}{=} Z \in \mathrm{RV}(\mu, \alpha)$. Let $\mathcal{G}_t = \sigma(\mathcal{A}, C_s, 0 \leq s \leq t)$ for $t \geq 0$, where $\mathcal{A}$ is independent of $\sigma(C_s, s \geq 0)$. Let $(H_t)_{t \geq 0}$ be predictable with respect to this filtration $d \times p$-matrix-valued process that is independent of $(C_t)$. Then the vector-valued stochastic integral

$$\int_0^T H_t \, \mathrm{d}C_t$$

has representation (1.1) with $A_j = H_{\tau_j} I\{1 \leq j \leq N_T\}$ and $Z_j = \Delta C_{\tau_j}$.

## 3. Convergence and tail behavior

Consider the random vector $X$ with stochastic representation (1.1). We will assume throughout most of this paper that

$$Z \in \mathrm{RV}(\mu, \alpha), \text{ and if } \alpha > 1, \text{ we assume additionally that } \mathrm{E}Z = 0. \qquad (3.1)$$

The assumption of the zero mean will allow us to work under relatively weak conditions. We will separately address what happens if the mean is not equal to zero.

For a matrix $A$, we denote by $\|A\|$ the operator norm of $A$. For a vector $z \in \mathbb{R}^d$, we denote the Euclidean norm by $|z|$.

We start by considering a linear process. That is, suppose $(A_j)_{j \in \mathbb{Z}}$ is a deterministic sequence of matrices. Then the following conditions are sufficient for a.s. convergence of the series (1.1):

$$\sum \|A_j\|^{\alpha - \varepsilon} < \infty, \qquad \text{for some } \varepsilon > 0, \qquad \text{if } 0 < \alpha \leq 2; \tag{3.2}$$

$$\sum \|A_j\|^2 < \infty, \qquad \qquad \text{if } \alpha > 2 \tag{3.3}$$

(here, and throughout the paper, we omit the summation index whenever it is clear what it is). This follows from Mikosch and Samorodnitsky (2000) in the case $d = p = 1$; the extension to the vector case is immediate. By Fubini's theorem, in the case of random $A_j$'s, *which are independent of the sequence* $(Z_j)$, the following conditions are, therefore, sufficient for a.s. convergence of the series (1.1):

$$\sum \|A_j\|^{\alpha - \varepsilon} < \infty \qquad \text{a.s. for some } \varepsilon > 0, \qquad \text{if } 0 < \alpha \leq 2; \tag{3.4}$$

$$\sum \|A_j\|^2 < \infty \qquad \text{a.s.} \qquad \qquad \text{if } \alpha > 2. \tag{3.5}$$

For linear processes, the conditions (3.2) and (3.3) turn out to be sufficient for the sum $X$ in (1.1) to acquire its regular variation from that of $(Z_j)$; in the case $d = p = 1$, this has been shown by Mikosch and Samorodnitsky (2000). It is clear that the conditions (3.4) and (3.5) will not suffice in the general case. If one looks at a general term in the sum (1.1), it is of the form $A_j Z_j$, and $Z_j$ is regularly varying. The tail behavior of such a product is usually controlled by a moment condition on the matrix $A_j$, of the type

$$\mathrm{E}[\|A_j\|^{\alpha + \varepsilon}] < \infty \quad \text{for some } \varepsilon > 0. \tag{3.6}$$

Then the term $A_j Z_j$ is regularly varying with limit measure $\mathrm{E}[\mu \circ A_j^{-1}(\cdot)]$ (see Basrak *et al.* (2002), Proposition A.1). This is the measure such that

$$\mathrm{E}[\mu \circ A_j^{-1}(B)] = \mathrm{E}[\mu\{z : A_j z \in B\}].$$

This result is usually referred to as *Breiman's theorem* (Breiman (1965)). Clearly, requiring (3.6) for each $j$ is too weak to control the tails of the infinite sum in (1.1). The above discussion shows, however, that we need to control both the small values of $A_j$'s that persist for long stretches of time, and the large values of $A_j$'s as well. This calls for a combination of different moment conditions. One such combination is presented in the following result.

**Theorem 3.1.** *Assume the 'predictable' framework (2.1)–(2.2). Suppose that (3.1) holds, $P(\bigcap_{j\geq 0}\{\|A_j\|=0\})=0$ and there is some $0<\varepsilon<\alpha$ such that*

$$\sum \mathrm{E}\|A_j\|^{\alpha-\varepsilon} < \infty \quad and \quad \sum \mathrm{E}\|A_j\|^{\alpha+\varepsilon} < \infty, \qquad if\ \alpha \in (0,1) \cup (1,2); \qquad (3.7)$$

$$\mathrm{E}\Big(\sum \|A_j\|^{\alpha-\varepsilon}\Big)^{(\alpha+\varepsilon)/(\alpha-\varepsilon)} < \infty, \qquad if\ \alpha \in \{1,2\}; \qquad (3.8)$$

$$\mathrm{E}\Big(\sum \|A_j\|^2\Big)^{(\alpha+\varepsilon)/2} < \infty, \qquad if\ \alpha \in (2,\infty). \qquad (3.9)$$

*Then the series (1.1) converges a.s. and*

$$\frac{P(u^{-1}X \in \cdot)}{P(|Z|>u)} \xrightarrow{v} \mathrm{E}\Big[\sum \mu \circ A_j^{-1}(\cdot)\Big], \qquad (3.10)$$

*as $u \to \infty$, on $\overline{\mathbb{R}}^d \setminus \{0\}$.*

**Remark 3.1.** Note that in the context of Theorem 3.1, it is only notationally different to consider a two-sided sum

$$X = \sum_{j \in \mathbb{Z}} A_j Z_j,$$

as long as one modifies the 'predictability' assumptions (2.1)–(2.2) appropriately. Indeed, the two-sided case can be turned into the one-sided case by relabelling the variables, for instance, by defining

$$A'_{2j} = A_j, \qquad Z'_{2j} = Z_j, \qquad j \geq 0,$$
$$A'_{-2j-1} = A_j, \qquad Z'_{-2j-1} = Z_j, \qquad j \leq -1.$$

For this relabelling, the predictability assumption is satisfied by assuming that $\sigma(Z_j, Z_{j+1}, \ldots)$ is independent of $\sigma(A_{-j}, \ldots, A_{-1}, A_0, \ldots, A_j)$ for $j \geq 0$ and that $\sigma(\ldots, Z_{j-1}, Z_j)$ is independent of $\sigma(A_j, \ldots, A_{-1}, A_0, \ldots, A_{-j-1})$ for $j \leq -1$. Indeed, in this case, one can take

$$\mathcal{F}_{2j} = \sigma(A_{-j}, \ldots, A_j, Z_{-j+2}, \ldots, Z_{j-1}),$$
$$\mathcal{F}_{2j+1} = \sigma(A_{-j-1}, \ldots, A_j, Z_{-j+1}, \ldots, Z_{j-1}).$$

If a different relabelling is used to go from the two-sided case to the one-sided, then one must modify the predictability assumption accordingly.

**Remark 3.2.** If one removes the assumption of zero mean in the case $\alpha > 1$, it is clear that the conclusion of Theorem 3.1 will still hold under the following additional assumption on the sequence $(A_j)$,

$$\text{the series } S_A = \sum A_j \text{ converges} \quad \text{and} \quad \lim_{u \to \infty} \frac{P(\|S_A\|>u)}{P(|Z|>u)} = 0, \qquad (3.11)$$

or under an even weaker assumption

$$\lim_{u\to\infty} \frac{P(|S_A \operatorname{E} Z| > u)}{P(|Z| > u)} = 0.$$

**Remark 3.3.** In the univariate case $p = d = 1$, the limiting measure $\mu$ of $Z$ has the representation

$$\mu(\mathrm{d}x) = w\alpha x^{-\alpha-1} I\{x > 0\}\,\mathrm{d}x + (1-w)\alpha |x|^{-\alpha-1} I\{x < 0\}\,\mathrm{d}x \qquad (3.12)$$

for some $w \in [0, 1]$. Then (3.10) becomes

$$\frac{P(X > u)}{P(|Z| > u)} \to \sum \operatorname{E}[|A_j|^\alpha (w I\{A_j > 0\} + (1-w) I\{A_j < 0\})].$$

**Remark 3.4.** The most general results on the tail behavior of the series (1.1) so far have been the works of Resnick and Willekens (1991) and Wang and Tang (2006) who considered the one-dimensional case $p = d = 1$ and the sequence $(A_j)$ being non-negative and independent of the sequence $(Z_j)$. Even in that particular case, our conditions in the case $\alpha \geq 1$ are strictly weaker (the conditions are identical for $0 < \alpha < 1$).

Next, we consider the implications of Theorem 3.1 in some of the examples presented in Section 1.

**Example 3.1 (Linear process).** Consider the linear process of Example 2.1. Since $(A_j)$ is a deterministic sequence, we immediately obtain the following statement.

**Corollary 3.1.** *Suppose that (3.1) holds and that there is some $0 < \varepsilon < \alpha$ such that*

$$\sum \|A_j\|^{\alpha-\varepsilon} < \infty, \qquad 0 < \alpha \leq 2,$$

$$\sum \|A_j\|^2 < \infty, \qquad \alpha > 2.$$

*Then the series (1.1) converges a.s. and*

$$\frac{P(u^{-1}X \in \cdot)}{P(|Z| > u)} \xrightarrow{v} \sum \mu \circ A_j^{-1}(\cdot)$$

*on $\overline{\mathbb{R}}^d \setminus \{0\}$.*

In the one-dimensional case $p = d = 1$ and $Z \in \mathrm{RV}(\alpha, \mu)$ with $\mu$ as in (3.12), we recover Mikosch and Samorodnitsky (2000), Lemma A.3 under identical assumptions.

**Example 3.2 (Random sum).** Consider the situation of Example 2.4: $d = p$ and $A_j = \mathrm{Id}\, I\{1 \leq j \leq N\}$, where $N$ is a positive-integer-valued random variable which is a stopping time with respect to a filtration $(\mathcal{G}_j, j \geq 0)$ to which the sequence $(Z_j)$ is adapted, and such that $\mathcal{G}_j$ is independent of $\sigma(Z_{j+1}, Z_{j+2}, \ldots)$ for $j \geq 0$. We then have the following result, which can be thought of as a 'tail Wald's identity'.

**Corollary 3.2.** *Suppose that (3.1) holds and that there is some $\tau > 0$ such that*

$$\mathrm{E}N < \infty \quad \text{if } \alpha \in (0,1) \cup (1,2), \tag{3.13}$$

$$\mathrm{E}N^{1+\tau} < \infty \quad \text{if } \alpha \in \{1,2\}, \tag{3.14}$$

$$\mathrm{E}N^{\alpha/2+\tau} < \infty \quad \text{if } \alpha \in (2,\infty). \tag{3.15}$$

*We then have*

$$\frac{P(u^{-1}X \in \cdot)}{P(|Z| > u)} \xrightarrow{v} \mathrm{E}N\mu(\cdot)$$

*on $\overline{\mathbb{R}}^d \setminus \{0\}$. If $\alpha > 1$, but $\mathrm{E}Z \neq 0$, then the same conclusion is obtained if one replaces (3.13), (3.14) and (3.15) by the assumption*

$$\lim_{u \to \infty} \frac{P(N > u)}{P(|Z| > u)} = 0. \tag{3.16}$$

In the one-dimensional case $d = p = 1$, with the noise variables $(Z_k)_{k \in \mathbb{Z}}$ being nonnegative, and independent of $N$, we recover the results of Stam (1973) and Fay *et al.* (2006).

***Example 3.3 (SRE).*** Consider the stationary solution to a stochastic recurrence equation of Example 2.2. From Theorem 3.1, we obtain the following result.

**Corollary 3.3.** *Suppose that $W \in \mathrm{RV}(\mu, \alpha)$ and that for some $\varepsilon > 0$,*

$$\mathrm{E}\|M\|^{\alpha+\varepsilon} < 1.$$

*Then the series (1.1) converges a.s. and (3.10) holds.*

**Proof.** For $\alpha \in (0,1) \cup (1,2)$, it follows by Jensen's inequality that

$$\mathrm{E}\sum \|A_j\|^{\alpha-\varepsilon} = \sum (\mathrm{E}\|M\|^{\alpha-\varepsilon})^j \leq \sum (\mathrm{E}\|M\|^{\alpha+\varepsilon})^{j(\alpha-\varepsilon)/(\alpha+\varepsilon)} < \infty.$$

Hence, (3.7) is satisfied. For $\alpha \in \{1,2\}$, it follows by convexity (see Kwapień and Woyczyński (1992), Lemma 3.3.1) that

$$\mathrm{E}\Big(\sum \|A_j\|^{\alpha-\varepsilon}\Big)^{(\alpha+\varepsilon)/(\alpha-\varepsilon)} \leq \Big(\sum (\mathrm{E}\|A_j\|^{\alpha+\varepsilon})^{(\alpha-\varepsilon)/(\alpha+\varepsilon)}\Big)^{(\alpha+\varepsilon)/(\alpha-\varepsilon)} \tag{3.17}$$

$$= \Big(\sum (\mathrm{E}\|M\|^{\alpha+\varepsilon})^{j(\alpha-\varepsilon)/(\alpha+\varepsilon)}\Big)^{(\alpha+\varepsilon)/(\alpha-\varepsilon)} < \infty. \tag{3.18}$$

Hence, (3.8) is satisfied. The case $\alpha > 2$ is similar. □

Note that we do not need to assume zero mean of the noise.
Specializing to the univariate case $p = d = 1$ and $\mu$ as in (3.12), we see that

$$\frac{P(X > x)}{P(|Z| > x)} \to \mathrm{E}\left[\sum_{j \geq 0} \mu[z : \Pi_{1-j,0} z > 1]\right]$$

$$= \mathrm{E}\left[\sum_{j \geq 0} \mu[z > 0 : (Y_{1-j} \cdots Y_0)^+ z > 1] + \mu[z < 0 : (Y_{1-j} \cdots Y_0)^- z < -1]\right]$$

$$= \sum_{j \geq 0} \mu[z > 1]\mathrm{E}[(Y_{1-j} \cdots Y_0)^+]^\alpha + \mu[z < -1]\mathrm{E}[(Y_{1-j} \cdots Y_0)^-]^\alpha$$

$$= \sum_{j \geq 0} w\mathrm{E}[(Y_{1-j} \cdots Y_0)^+]^\alpha + (1-w)\mathrm{E}[(Y_{1-j} \cdots Y_0)^-]^\alpha$$

$$= w \sum_{j \geq 0} \sum_{k=0,1,\ldots,j,\text{even}} \binom{j}{k} (\mathrm{E}(Y^-)^\alpha)^k (\mathrm{E}(Y^+)^\alpha)^{j-k}$$

$$+ (1-w) \sum_{j \geq 0} \sum_{k=0,1,\ldots,j,\text{ odd}} \binom{j}{k} (\mathrm{E}(Y^-)^\alpha)^k (\mathrm{E}(Y^+)^\alpha)^{j-k}$$

$$= w \sum_{k \text{ even}} \frac{(\mathrm{E}(Y^-)^\alpha)^k}{(1 - \mathrm{E}(Y^+)^\alpha)^{k+1}} + (1-w) \sum_{k \text{ odd}} \frac{(\mathrm{E}(Y^-)^\alpha)^k}{(1 - \mathrm{E}(Y^+)^\alpha)^{k+1}}$$

$$= \frac{w(1 - \mathrm{E}(Y^+)^\alpha) + (1-w)\mathrm{E}(Y^-)^\alpha}{(1 - \mathrm{E}(Y^+)^\alpha)^2 - (\mathrm{E}(Y^-)^\alpha)^2}.$$

In particular, if $Y$ is non-negative, then this reduces to

$$\frac{P(X > x)}{P(|Z| > x)} \to w(1 - \mathrm{E}Y^\alpha)^{-1}.$$

In this particular case, we recover (under identical assumptions) the results of Grey (1994) or Konstantinides and Mikosch (2005), Proposition 2.2.

## 4. Proof of Theorem 3.1

The outline of the proof is as follows. First, we show that the conditions (3.7)–(3.9) guarantee almost sure convergence of the sum in (1.1). Then, to prove (3.10), we start with the univariate case and divide into four cases: $0 < \alpha < 1$ in Section 4.1; $\alpha = 1$ in Section 4.2; $1 < \alpha < 2$ in Section 4.3; $\alpha \geq 2$ in Section 4.4. Finally, the multidimensional case is treated in Section 4.5.

We start with the almost sure convergence of the sum (1.1). Let $(\hat{Z}_j)$ be a sequence with the same law as $(Z_j)$, independent of the sequence $(A_j)$ (we may need to enlarge the underlying probability space to construct such a sequence). The series $\sum_{j=0}^{\infty} A_j \hat{Z}_j$ converges a.s. by Fubini's theorem and (3.4)–(3.5). Furthermore, the sequences $(A_j \hat{Z}_j)$ and $(A_j Z_j)$ are tangent with respect to the filtration $\tilde{\mathcal{F}}_j = \sigma(\mathcal{F}_{j+1}, \hat{Z}_0, \ldots, \hat{Z}_j)$ for $j \geq 0$. That is, $\text{Law}(A_j Z_j \mid \tilde{\mathcal{F}}_{j-1}) = \text{Law}(A_j \hat{Z}_j \mid \tilde{\mathcal{F}}_{j-1})$. Also, note that the sequence $(A_j \hat{Z}_j)$ is conditionally independent, given the $\sigma$-field $\mathcal{G} = \sigma((A_k)_{k \geq 0}, (Z_k)_{k \geq 0})$ (see Kwapień and Woyczyński (1992), Section 4.3 for details). Therefore, Kwapień and Woyczyński (1992), Corollary 5.7.1 guarantees a.s. convergence of the series in (1.1).

We will prove now (3.10). We start with the univariate case $d = p = 1$, in which case the statement of the proposition reduces to

$$\lim_{x \to \infty} \frac{P(X > x)}{P(|Z| > x)} = \sum_{j=1}^{\infty} (w\mathrm{E}((A_j)^+)^\alpha + (1-w)\mathrm{E}((A_j)^-)^\alpha) \tag{4.1}$$

with $w = \lim_{x \to \infty} P(Z > x)/P(|Z| > x)$ and where $a^+$ and $a^-$ are, respectively, the positive part and the negative part of a real number $a$.

For finite $n \geq 1$ as $x \to \infty$, we have

$$\frac{P(\sum_{j=1}^{n} A_j Z_j > x)}{P(|Z| > x)} \to \sum_{j=1}^{n} (w\mathrm{E}((A_j)^+)^\alpha + (1-w)\mathrm{E}((A_j)^-)^\alpha),$$

by Lemma 4.3 below. Therefore, it is sufficient to show that

$$\lim_{n \to \infty} \limsup_{x \to \infty} \frac{P(\sum_{j > n} A_j Z_j > x)}{P(|Z| > x)} = 0. \tag{4.2}$$

### 4.1. The case $0 < \alpha < 1$

We start with the case $0 < \alpha < 1$, in which case we will actually check that (4.2) holds with the sum of absolute values in the numerator. The first step is to show that for any $M > 0$, (4.2) holds with each $A_j$ replaced by $\tilde{A}_j = A_j I\{|A_j| < M\}$ and, by scaling, it is enough to consider the case $M = 1$. We may decompose the probability as

$$\frac{P(\sum |\tilde{A}_j Z_j| > x)}{P(|Z| > x)} = \underbrace{\frac{P(\sum |\tilde{A}_j Z_j| > x, \bigvee |\tilde{A}_j Z_j| > x)}{P(|Z| > x)}}_{\text{I}} + \underbrace{\frac{P(\sum |\tilde{A}_j Z_j| > x, \bigvee |\tilde{A}_j Z_j| \leq x)}{P(|Z| > x)}}_{\text{II}}, \tag{4.3}$$

with $\bigvee$ denoting maximum. For I, an upper bound can be constructed as

$$\frac{P(\sum|\tilde{A}_j Z_j| > x, \bigvee|\tilde{A}_j Z_j| > x)}{P(|Z| > x)} \leq \frac{P(\bigvee|\tilde{A}_j Z_j| > x)}{P(|Z| > x)}$$

$$\leq \sum \frac{P(|\tilde{A}_j Z_j| > x)}{P(|Z| > x)}$$

$$= \sum \frac{\int_0^1 P(y|Z_j| > x) P(|\tilde{A}_j| \in \mathrm{d}y)}{P(|Z| > x)}.$$

Using Potter's bounds (see, e.g., Resnick (1987)) it follows that there exists $x_0 > 0$ such that $P(|Z| > x/y)/P(|Z| > x) \leq cy^{\alpha - \varepsilon}$ for $x > x_0$ and $0 < y \leq 1$. Hence, the last expression is bounded above by

$$c \sum_{j>n} \int_0^1 y^{\alpha - \varepsilon} P(|A_j| \in \mathrm{d}y) \leq c \sum_{j>n} \mathrm{E}|A_j|^{\alpha - \varepsilon} \to 0 \qquad \text{as } n \to \infty,$$

by (3.7). For II, Markov's inequality implies

$$\mathrm{II} \leq \frac{P(\sum|\tilde{A}_j Z_j| I\{|\tilde{A}_j Z_j| \leq x\} > x)}{P(|Z| > x)}$$

$$\leq \sum \frac{\mathrm{E}[|\tilde{A}_j Z_j| I\{|\tilde{A}_j Z_j| \leq x\}]}{x P(|Z| > x)}$$

$$= \sum \int_0^1 y \frac{\mathrm{E}[|Z_j| I\{|Z_j| \leq x/y\}]}{x P(|Z| > x)} P(|A_j| \in \mathrm{d}y).$$

By Karamata's theorem (see Resnick (1987))

$$\mathrm{E}[|Z| I\{|Z| \leq x\}] \sim \alpha (1-\alpha)^{-1} x P(|Z| > x)$$

and there exists $x_0$ such that for $x > x_0$, the last expression is bounded from above by

$$c \sum_{j>n} \int_0^1 y y^{-1+\alpha-\varepsilon} P(|A_j| \in \mathrm{d}y)$$

$$\leq c \sum_{j>n} \mathrm{E}[|A_j|^{\alpha-\varepsilon}] \to 0 \qquad \text{as } n \to \infty,$$

by (3.7). Combining I and II proves (4.2) for $(\tilde{A}_j)$.

Next, for $M > 0$ and $v < x/M$, the remaining term can be bounded by

$$P\Big(\sum |A_j I\{|A_j| > M\} Z_j| > x\Big)$$

$$\leq P\Big(\sum |A_j| I\{M < |A_j| < x/v\} |Z_j| > x/2\Big)$$
$$+ P\Big(\sum |A_j| I\{|A_j| \geq x/v\} |Z_j| > x/2\Big).$$

The result will follow once we show that, uniformly in $n$,

$$\lim_{M \to \infty} \limsup_{x \to \infty} \frac{P(\sum_{j>n} |A_j| I\{M \leq |A_j| < x/v\} |Z_j| > x)}{P(|Z| > x)} = 0 \quad (4.4)$$

and

$$\lim_{M \to \infty} \limsup_{x \to \infty} \frac{P(\sum_{j>n} |A_j| I\{|A_j| \geq x/v\} |Z_j| > x)}{P(|Z| > x)} = 0. \quad (4.5)$$

Let us start with (4.4). Using the decomposition which lead to I and II above, we see that it is sufficient to show that

$$\lim_{M \to \infty} \limsup_{x \to \infty} \sum \frac{P(|A_j| I\{M \leq |A_j| < x/v\} |Z_j| > x)}{P(|Z| > x)} = 0 \quad (4.6)$$

and

$$\lim_{M \to \infty} \limsup_{x \to \infty} \frac{P(\sum |A_j| I\{M \leq |A_j| < x/v\} |Z_j| I\{|A_j Z_j| \leq x\} > x)}{P(|Z| > x)} = 0. \quad (4.7)$$

The term in (4.6) can be written as

$$\sum \frac{P(|A_j| I\{M \leq |A_j| < x/v\} |Z_j| > x)}{P(|Z| > x)}$$
$$= \sum \int_M^{x/v} \frac{P(|Z_j| > x/y)}{P(|Z| > x)} P(|A_j| \in \mathrm{d}y).$$

For $v$ and $x$ sufficiently large, Potter's bounds imply that this is bounded above by

$$c \sum \int_M^{x/v} y^{\alpha+\varepsilon} P(|A_j| \in \mathrm{d}y) \leq c \sum \mathrm{E}|A_j|^{\alpha+\varepsilon} I\{A_j \geq M\}.$$

Since $\sum \mathrm{E}|A_j|^{\alpha+\varepsilon} < \infty$, by assumption (3.7), the sum converges to zero as $M \to \infty$.

We now turn to (4.7). Markov's inequality gives

$$\frac{P(\sum |A_j| I\{M \leq |A_j| < x/v\} |Z_j| I\{|A_j Z_j| \leq x\} > x)}{P(|Z| > x)}$$
$$\leq \sum \frac{\mathrm{E}[|A_j| I\{M \leq |A_j| < x/v\} |Z_j| I\{|A_j Z_j| \leq x\}]}{x P(|Z| > x)}$$

$$= \sum \int_M^{x/v} \frac{y\mathrm{E}[|Z_j|I\{|Z_j| \leq x/y\}]}{xP(|Z| > x)} P(|A_j| \in \mathrm{d}y).$$

Lemma 4.1 implies that for $v$ and $x$ sufficiently large, we have the upper bound

$$c\sum \int_M^{x/v} y^{\alpha+\varepsilon} P(|A_j| \in \mathrm{d}y) \leq c\sum \mathrm{E}|A_j|^{\alpha+\varepsilon} I\{|A_j| \geq M\},$$

which converges to zero.

Finally, we want to show (4.5). By Potter's bounds, $P(|Z| > x) \geq x^{-\alpha-\varepsilon}$ for $x$ sufficiently large. Hence,

$$\frac{P(\sum |A_j|I\{|A_j| \geq x/v\}|Z_j| > x)}{P(|Z| > x)} \leq \frac{P(\bigvee |A_j| > x/v)}{P(|Z| > x)}$$

$$\leq \sum \frac{P(|A_j| > x/v)}{P(|Z| > x)}$$

$$\leq \sum x^{\alpha+\varepsilon} P(|A_j| > x/v)$$

$$\leq v^{\alpha+\varepsilon} \sum \int_{x/v}^\infty y^{\alpha+\varepsilon} P(|A_j| \in \mathrm{d}y)$$

$$\leq v^{\alpha+\varepsilon} \sum \mathrm{E}|A_j|^{\alpha+\varepsilon} I\{|A_j| \geq x/v\},$$

which converges to zero as $x \to \infty$. This completes the proof of (4.2) in the case $0 < \alpha < 1$.

### 4.2. The case $\alpha = 1$

As in the case $0 < \alpha < 1$, we start by proving (4.2) with each $A_j$ replaced by $\tilde{A}_j = A_j I\{|A_j| < M\}$. Again, by scaling, it is enough to consider the case $M = 1$. We use the decomposition (4.3). The argument for I is the same as for $0 < \alpha < 1$, so it is sufficient to consider II. Let

$$L(y) = \mathrm{E}|Z|I\{|Z| \leq y\}, \qquad y > 0.$$

Note that since $\alpha = 1$, $L$ is slowly varying at infinity. Write

$$P\Big(\sum |\tilde{A}_j Z_j| I\{|\tilde{A}_j Z_j| \leq x\} > x\Big)$$

$$\leq P\Big(\sum |\tilde{A}_j|[|Z_j|I\{|\tilde{A}_j Z_j| \leq x\} - L(x/|\tilde{A}_j|)] > x/2\Big)$$

$$+ P\Big(\sum |\tilde{A}_j|L(x/|\tilde{A}_j|) > x/2\Big) := \mathrm{II}_a + \mathrm{II}_b.$$

Consider first $\mathrm{II}_a$. Let $0 < \delta < 1$. By Markov's inequality,

$$P\Big(\sum |\tilde{A}_j|[|Z_j|I\{|\tilde{A}_j Z_j| \leq x\} - L(x/|\tilde{A}_j|)] > x/2\Big)$$

$$\leq \left(\frac{2}{x}\right)^{1+\delta} \mathrm{E}\Big|\sum |\tilde{A}_j|[|Z_j|I\{|\tilde{A}_j Z_j| \leq x\} - L(x/|\tilde{A}_j|)]\Big|^{1+\delta}.$$

By predictability of $(\tilde{A}_j)$, the above sum is a sum of martingale differences. By the Burkholder–Davis–Gundy inequality (see, e.g., Protter (2004)), the last expression is bounded above by

$$\frac{c}{x^{1+\delta}} \mathrm{E}\Big|\sum |\tilde{A}_j|^2 [|Z_j|I\{|\tilde{A}_j Z_j| \leq x\} - L(x/|\tilde{A}_j|)]^2\Big|^{(1+\delta)/2}$$
$$\leq \frac{c}{x^{1+\delta}} \sum \mathrm{E}||\tilde{A}_j|[|Z_j|I\{|\tilde{A}_j Z_j| \leq x\} - L(x/|\tilde{A}_j|)]|^{(1+\delta)}.$$

Conditioning on the $A_j$ variables, we see that

$$\frac{\mathrm{II}_a}{P(|Z|>x)} \leq c \sum \int_0^1 y^{1+\delta} \frac{\mathrm{E}||Z_j|I\{|Z_j| \leq x/y\} - L(x/y)|^{1+\delta}}{x^{1+\delta} P(|Z|>x)} P(|A_j| \in \mathrm{d}y)$$
$$\leq c \sum \int_0^1 y^{1+\delta} \frac{\mathrm{E}|Z_j|^{1+\delta} I\{|Z_j| \leq x/y\} + L(x/y)^{1+\delta}}{x^{1+\delta} P(|Z|>x)} P(|A_j| \in \mathrm{d}y)$$
$$\leq 2c \sum \int_0^1 y^{1+\delta} \frac{\mathrm{E}|Z_j|^{1+\delta} I\{|Z_j| \leq x/y\}}{x^{1+\delta} P(|Z|>x)} P(|A_j| \in \mathrm{d}y).$$

Once again, by Karamata's theorem (Resnick (1987)),

$$\mathrm{E}[|Z|^{1+\delta} I\{|Z| \leq x\}] \sim \alpha \delta^{-1} x^{1+\delta} P(|Z|>x) \tag{4.8}$$

and there exists $x_0$ such that for $x > x_0$, the above expression is bounded by

$$c \sum_{j>n} \int_0^1 y^{1+\delta} y^{-(1+\delta)+1-\varepsilon} P(|A_j| \in \mathrm{d}y) \leq c \sum_{j>n} \mathrm{E}[|A_j|^{1-\varepsilon}] \to 0 \quad \text{as } n \to \infty,$$

by (3.8). This takes care of the term $\mathrm{II}_a$. For $\mathrm{II}_b$, we have, for $\varepsilon > 0$, by slow variation of $L(x)$, using Markov's inequality and (3.8),

$$\frac{\mathrm{II}_b}{P(|Z|>x)} \leq \frac{P(\sum |\tilde{A}_j|(x/|\tilde{A}_j|)^\varepsilon > cx)}{P(|Z|>x)}$$
$$= \frac{P(\sum |\tilde{A}_j|^{1-\varepsilon} > cx^{1-\varepsilon})}{P(|Z|>x)}$$
$$\leq \frac{C}{x^{(1+\varepsilon)} P(|Z|>x)} \mathrm{E}\Big(\sum |\tilde{A}_j|^{1-\varepsilon}\Big)^{(1+\varepsilon)/(1-\varepsilon)} \to 0 \quad \text{as } x \to \infty.$$

Therefore, we have (4.2) for $(\tilde{A}_j)$.

Further, as in the case $0 < \alpha < 1$, we need to check (4.4) and (4.5). The argument for (4.5) works without changes in the present case and the same is true for the first half of the argument for (4.4), presented in (4.6). Therefore, it remains only to consider the second half of (4.4), namely to prove that for $v$ large enough,

$$\lim_{M\to\infty}\limsup_{x\to\infty}\frac{P(\sum|A_j|I\{M\leq|A_j|<x/v\}|Z_j|I\{|Z_j|\leq x/|A_j|\}>x)}{P(|Z|>x)}=0.$$

The argument is similar to the one used above. First, we can decompose the probability

$$P\Big(\sum|A_j|I\{M\leq|A_j|<x/v\}|Z_j|I\{|Z_j|\leq x/|A_j|\}>x\Big)$$
$$\leq P\Big(\sum|A_j|I\{M\leq|A_j|<x/v\}(|Z_j|I\{|Z_j|\leq x/|A_j|\}-L(x/|A_j|))>x/2\Big)$$
$$+P\Big(\sum|A_j|I\{M\leq|A_j|<x/v\}L(x/|A_j|)>x/2\Big):=\mathrm{III}_a+\mathrm{III}_b.$$

Using Markov's inequality and the Burkholder–Davis–Gundy inequality we see that for $0<\delta<1$,

$$\mathrm{III}_a\leq\frac{c}{x^{1+\delta}}\sum\mathrm{E}|A_j|^{1+\delta}I\{M\leq|A_j|<x/v\}||Z_j|I\{|Z_j|\leq x/|A_j|\}-L(x/|A_j|)|^{1+\delta}.$$

Notice that by (4.8) and for $v$ large enough,

$$\mathrm{E}|A_j|^{1+\delta}I\{M\leq|A_j|<x/v\}||Z_j|I\{|Z_j|\leq x/A_j\}-L(x/|A_j|)|^{1+\delta}$$
$$=\mathrm{E}(|A_j|^{1+\delta}I\{M\leq|A_j|<x/v\}\mathrm{E}_Z||Z_j|I\{|Z_j|\leq x/|A_j|\}-L(x/|A_j|)|^{1+\delta})$$
$$\leq c\mathrm{E}\bigg(|A_j|^{1+\delta}I\{M\leq|A_j|<x/v\}\bigg(\frac{x}{|A_j|}\bigg)^{1+\delta}P_Z(|Z|>x/|A_j|)\bigg)$$
$$\leq cx^{1+\delta}P(|Z|>x)\mathrm{E}(|A_j|^{1+\varepsilon}I\{|A_j|\geq M\})$$

for some $\varepsilon>0$. Here, $P_Z$ and $\mathrm{E}_Z$ indicate that the probability and expectation, respectively, are computed with respect to the $Z$ variables (i.e., conditionally on $A_j$). In the last step, we used independence and the fact that $P(|Z|>z)$ is regularly varying. We now see that

$$\frac{\mathrm{III}_a}{P(|Z|>x)}\leq c\sum\mathrm{E}|A_j|^{1+\varepsilon}I\{|A_j|\geq M\},$$

which converges to zero by (3.8).

It remains to consider $\mathrm{III}_b$. For $\mathrm{III}_b$, we have, for $\varepsilon>0$, by slow variation of $L(x)$, using Markov's inequality and (3.8),

$$\frac{\mathrm{III}_b}{P(|Z|>x)}\leq\frac{P(\sum|A_j|I\{M\leq|A_j|<x/v\}(x/|A_j|)^\varepsilon>cx)}{P(|Z|>x)}$$
$$=\frac{P(\sum|A_j|^{1-\varepsilon}I\{M\leq|A_j|<x/v\}>cx^{1-\varepsilon})}{P(|Z|>x)}$$

$$\leq \frac{C}{x^{(1+\varepsilon)}P(|Z|>x)}\mathrm{E}\Big(\sum|A_j|^{1-\varepsilon}\Big)^{(1+\varepsilon)/(1-\varepsilon)} \to 0 \quad \text{as } x\to\infty.$$

This completes the proof of (4.2) for $\alpha=1$.

### 4.3. The case $1<\alpha<2$

Next, let $1<\alpha<2$. Assume first that the law of $Z$ is continuous and that the limit measure $\mu$ assigns positive weights to both $(-\infty,0)$ and $(0,\infty)$.

We start, as above, in the case $0<\alpha<1$ (this time, the absolute values stay outside the sum) and establish (4.2) for $\tilde{A}_j$. We start with the same decomposition as in (4.3), except that we decompose not according to whether or not $\bigvee|\tilde{A}_j Z_j|>x$, but rather whether or not for some $j$, $Z_j$ does not belong to the interval $[-h(x/\tilde{A}_j), x/\tilde{A}_j]$ for the function $h$ in Lemma 4.2. The idea here is that $h$ is constructed such that $ZI\{-h(y)\leq Z \leq y\}$ has mean zero so that martingale inequalities can be applied. Obviously, the argument for I in the decomposition (4.3) works for any $\alpha>1$ and so it covers the case $1<\alpha<2$. Let us consider the term II in that decomposition. Let $0<\delta<2-\alpha$. By Markov's inequality,

$$\begin{aligned}
&P\Big(\Big|\sum \tilde{A}_j Z_j I\{-h(x/\tilde{A}_j)\leq Z_j \leq x/\tilde{A}_j\}\Big|>x\Big) \\
&\leq \frac{1}{x^{\alpha+\delta}}\mathrm{E}\Big|\sum \tilde{A}_j Z_j I\{-h(x/\tilde{A}_j)\leq Z_j \leq x/\tilde{A}_j\}\Big|^{\alpha+\delta}.
\end{aligned} \quad (4.9)$$

By the predictability, zero mean of the $Z$'s and the property of the function $h$, the above is a sum of martingale differences. Therefore, we can use the Burkholder–Davis–Gundy inequality (see, e.g., Protter (2004)) to conclude that for large $x$, the above is bounded by

$$\frac{c}{x^{\alpha+\delta}}\mathrm{E}\Big|\sum (\tilde{A}_j Z_j I\{-h(x/\tilde{A}_j)\leq Z_j \leq x/\tilde{A}_j\})^2\Big|^{(\alpha+\delta)/2}$$
$$\leq \frac{c}{x^{\alpha+\delta}}\sum \mathrm{E}|\tilde{A}_j Z_j I\{-h(x/\tilde{A}_j)\leq Z_j \leq x/\tilde{A}_j\}|^{\alpha+\delta},$$

where $c$ is a finite positive constant that is allowed to change in the sequel. In the last step, we used the fact that the $\ell_q$-norm is bounded by the $\ell_p$-norm for $p<q$, that is, $\|\cdot\|_{\ell_q}\leq\|\cdot\|_{\ell_p}$. Conditioning on the $A_j$ variables, we see that

$$\frac{\text{II}}{P(|Z|>x)}\leq c\sum\int_0^1 y^{\alpha+\delta}\frac{\mathrm{E}[|Z_j|^{\alpha+\delta}I\{-h(x/y)\leq Z_j \leq x/y\}]}{x^{\alpha+\delta}P(|Z|>x)}P(|A_j|\in\mathrm{d}y).$$

Once again, by Karamata's theorem (Resnick (1987)),

$$\mathrm{E}[|Z|^{\alpha+\delta}I\{|Z|\leq x\}] \sim \alpha\delta^{-1}x^{\alpha+\delta}P(|Z|>x) \qquad (4.10)$$

and there exists some $x_0$ such that for $x > x_0$, the above expression is bounded by

$$c \sum_{j>n} \int_0^1 y^{\alpha+\delta} y^{-(\alpha+\delta)+\alpha-\varepsilon} P(|A_j| \in \mathrm{d}y) \leq c \sum_{j>n} \mathrm{E}[|A_j|^{\alpha-\varepsilon}] \to 0 \quad \text{as } n \to \infty,$$

by (3.7). Therefore, we have (4.2) for $(\tilde{A}_j)$.

Further, as in the case $0 < \alpha < 1$, we need to check (4.4) and (4.5). The argument for (4.5) works without changes in the present case and the same is true for the first half of the argument for (4.4), presented in (4.6). Therefore, it remains only to consider the second half of (4.4), namely to prove that for $v$ large enough,

$$\lim_{M \to \infty} \limsup_{x \to \infty} \frac{P(|\sum A_j I\{M \leq |A_j| < x/v\} Z_j I\{-h(x/A_j) \leq Z_j \leq x/A_j\}| > x)}{P(|Z| > x)} = 0.$$

The argument is similar to the one used above. Using Markov's inequality, predictability, zero mean and the definition of the function $h$ allows us, once again, to use the Burkholder–Davis–Gundy inequality and see that for $0 < \delta < 2 - \alpha$,

$$P\left(\left|\sum A_j I\{M \leq |A_j| < x/v\} Z_j I\{-h(x/A_j) \leq Z_j \leq x/A_j\}\right| > x\right)$$
$$\leq \frac{c}{x^{\alpha+\delta}} \sum \mathrm{E}|A_j Z_j I\{M \leq |A_j| < x/v\} I\{-h(x/A_j) \leq Z_j \leq x/A_j\}|^{\alpha+\delta} \quad (4.11)$$
$$:= \mathrm{III}.$$

Notice that by (4.10), for and $v$ large enough and using the fact that $h(x) \leq cx$,

$$\mathrm{E}|A_j Z_j I\{M \leq |A_j| < x/v\} I\{-h(x/A_j) \leq Z_j \leq x/A_j\}|^{\alpha+\delta}$$
$$= \mathrm{E}(|A_j|^{\alpha+\delta} I\{M \leq |A_j| < x/v\} \mathrm{E}_Z |Z_j|^{\alpha+\delta} I\{-h(x/A_j) \leq Z_j \leq x/A_j\})$$
$$\leq c\mathrm{E}\left(|A_j|^{\alpha+\delta} I\{M \leq |A_j| < x/v\} \left(\frac{x}{|A_j|}\right)^{\alpha+\delta} P_Z(|Z| > x/|A_j|)\right)$$
$$\leq cx^{\alpha+\delta} P(|Z| > x) \mathrm{E}(|A_j|^{\alpha+\varepsilon} I\{|A_j| \geq M\}).$$

In the last step, we used independence and the fact that $P(|Z| > z)$ is regularly varying. We now see that

$$\frac{\mathrm{III}}{P(|Z| > x)} \leq c \sum \mathrm{E}|A_j|^{\alpha+\varepsilon} I\{|A_j| \geq M\},$$

which converges to zero by (3.7). Therefore, we have proven (4.2) for $1 < \alpha < 2$, under the additional assumption that the law of $Z$ is continuous and that the limit measure $\mu$ assigns positive weights to both $(-\infty, 0)$ and $(0, \infty)$. In general, let $(\tilde{Z}_j)$ be an i.i.d. sequence, independent of the sequences $(A_j)$ and $(Z_j)$ (we enlarge the probability space if necessary), such that each $\tilde{Z}_j$ is continuous, symmetric and

$$\lim_{x \to \infty} \frac{P(|\tilde{Z}_j| > x)}{P(|Z| > x)} = 1.$$

By symmetry and independence,

$$P\Big(\Big|\sum A_j Z_j\Big| > x\Big) \leq 2P\Big(\Big|\sum A_j(Z_j + \tilde{Z}_j)\Big| > x\Big).$$

However, the sequence $(Z_j + \tilde{Z}_j)$ satisfies the extra assumptions and so we have established (4.2) for $1 < \alpha < 2$ in full generality.

### 4.4. The case $\alpha \geq 2$

The proof of the relation (4.2) for $\alpha \geq 2$ proceeds similarly to the case $1 < \alpha < 2$. Specifically, we need to estimate both the term II in (4.9) and the term III in (4.11). As in the case $1 < \alpha < 2$, we may, and will, assume that the random variables $(Z_j)$ satisfy the assumptions of Lemma 4.2. Furthermore, using Kwapień and Woyczyński (1992), part (iv) of Theorem 5.2.1, we may assume that the sequence $(Z_j)$ is independent of the sequence $(A_j)$. Indeed, to achieve that, we simply replace the sequence $(Z_j)$ with the sequence $(\hat{Z}_j)$ defined at the beginning of the proof of the theorem and use the tangency.

We start with the case $\alpha = 2$. We first estimate II in (4.9). Starting with the Burkholder–Davis–Gundy inequality as before, we proceed as follows. Using Jensen's inequality,

$$\sum(\tilde{A}_j Z_j I\{-h(x/\tilde{A}_j) \leq Z_j \leq x/\tilde{A}_j\})^2$$
$$\leq \Big(\sum \tilde{A}_j^2\Big)^{\delta/(2+\delta)} \Big(\sum \tilde{A}_j^2 |Z_j|^{2+\delta} I\{-h(x/\tilde{A}_j) \leq Z_j \leq x/\tilde{A}_j\}\Big)^{2/(2+\delta)}$$

and so, by conditioning,

$$\text{II} \leq \frac{c}{x^{2+\delta}} \text{E}\Big[\Big(\sum \tilde{A}_j^2\Big)^{\delta/2} \sum \tilde{A}_j^2 |Z_j|^{2+\delta} I\{-h(x/\tilde{A}_j) \leq Z_j \leq x/\tilde{A}_j\}\Big]$$
$$= \frac{c}{x^{2+\delta}} \text{E}\Big[\Big(\sum \tilde{A}_j^2\Big)^{\delta/2} \sum \tilde{A}_j^2 \text{E}_Z(|Z_j|^{2+\delta} I\{-h(x/\tilde{A}_j) \leq Z_j \leq x/\tilde{A}_j\})\Big].$$

Using (4.8), we see that for large $x$, this expression is bounded from above by

$$\frac{c}{x^{2+\delta}} \text{E}\Big[\Big(\sum \tilde{A}_j^2\Big)^{\delta/2} \sum \tilde{A}_j^2 \Big(\frac{x}{|\tilde{A}_j|}\Big)^{2+\delta} P(|Z| > x/|\tilde{A}_j|)\Big]$$
$$= c\text{E}\Big[\Big(\sum \tilde{A}_j^2\Big)^{\delta/2} \sum |\tilde{A}_j|^{-\delta} P(|Z| > x/|\tilde{A}_j|)\Big]$$
$$\leq cP(|Z| > x)\text{E}\Big[\Big(\sum \tilde{A}_j^2\Big)^{\delta/2} \sum |\tilde{A}_j|^{2-\varepsilon-\delta}\Big]$$
$$\leq cP(|Z| > x)\text{E}\Big[\Big(\sum |\tilde{A}_j|^{2-\varepsilon-\delta}\Big)^{1+\delta/(2-\varepsilon-\delta)}\Big]$$

(when writing sums as above, and in the sequel, we adopt the convention of not including the terms with $\tilde{A}_j = 0$). In the last two steps, we used regular variation of $P(|Z| > z)$ and $\|\cdot\|_{\ell_q} \leq \|\cdot\|_{\ell_p}$, for $p < q$, respectively. Choosing $\delta$ and $\varepsilon$ small enough and using (3.8), we see that

$$\lim_{n \to \infty} \limsup_{x \to \infty} \frac{\mathrm{II}}{P(|Z| > x)} = 0.$$

The argument for III is similar; we present the main steps. Write

$$\mathrm{III} \leq \frac{c}{x^{\alpha+\delta}} \mathrm{E}\Big[\Big(\sum A_j^2 I\{|A_j| \geq M\}\Big)^{\delta/2}$$
$$\times \sum A_j^2 I\{M \leq |A_j| < x/v\} |Z_j|^{2+\delta} I\{-h(x/A_j) \leq Z_j \leq x/A_j\}\Big]$$
$$\leq c\mathrm{E}\Big[\Big(\sum A_j^2 I\{|A_j| \geq M\}\Big)^{\delta/2}$$
$$\times \sum |A_j|^{-\delta} I\{M \leq |A_j| < x/v\} P(|Z| > x/|A_j|)\Big]$$
$$\leq cP(|Z| > x) \mathrm{E}\Big[\Big(\sum |A_j|^{2+\varepsilon-\delta} I\{|A_j| \geq M\}\Big)^{1+\delta/(2+\varepsilon-\delta)}\Big].$$

Choosing, for example, $\delta = \varepsilon$ small enough, we see by (3.8) that

$$\lim_{n \to \infty} \limsup_{x \to \infty} \frac{\mathrm{III}}{P(|Z| > x)} = 0.$$

This establishes the statement (4.2) for $\alpha = 2$.

Next, we look at the case $\alpha > 2$ and $\alpha$ not equal to an even integer. We start with the term II. Let $k = \lceil(\alpha + \delta)/2\rceil$, that is, the smallest integer greater than or equal to $(\alpha + \delta)/2$. Our assumption on $\alpha$ implies that for $\delta$ small enough,

$$2(k-1) < \alpha. \tag{4.12}$$

Proceeding from the Burkholder–Davis–Gundy bound on II, we have

$$\mathrm{II} \leq \frac{c}{x^{\alpha+\delta}} \mathrm{E}\Bigg[\sum_{j_1} \cdots \sum_{j_k} \tilde{A}_{j_1}^2 \cdots \tilde{A}_{j_k}^2 Z_{j_1}^2 \cdots Z_{j_k}^2$$
$$\times I\{-h(x/\tilde{A}_{j_i}) \leq Z_{j_i} \leq x/\tilde{A}_{j_i} \text{ for } i = 1, \ldots, k\}\Bigg]^{(\alpha+\delta)/(2k)}$$
$$\leq \frac{c}{x^{\alpha+\delta}} \mathrm{E}\Big[\sum |\tilde{A}_j|^{2k} |Z_j|^{2k} I\{-h(x/A_j) \leq Z_j \leq x/A_j\}\Big]^{(\alpha+\delta)/(2k)}$$
$$+ \frac{c}{x^{\alpha+\delta}} \mathrm{E}\Bigg[\sum_{(j_1,\ldots,j_k) \in D_k} \tilde{A}_{j_1}^2 \cdots \tilde{A}_{j_k}^2 Z_{j_1}^2 \cdots Z_{j_k}^2$$

$$\times I\{-h(x/\tilde{A}_{j_i}) \leq Z_{j_i} \leq x/\tilde{A}_{j_i} \text{ for } i=1,\ldots,k\}\Bigg]^{(\alpha+\delta)/(2k)}$$

$$:= \mathrm{II}_a + \mathrm{II}_b,$$

where $D_k = \{(j_1,\ldots,j_k) \text{ such that not all } j_1,\ldots,j_k \text{ are equal}\}$. Note that by the definition of $k$ and, again, $\|\cdot\|_{\ell_q} \leq \|\cdot\|_{\ell_p}$, for $p < q$,

$$\mathrm{II}_a \leq \frac{c}{x^{\alpha+\delta}} \mathrm{E} \sum |\tilde{A}_j|^{\alpha+\delta} |Z_j|^{\alpha+\delta} I\{-h(x/A_j) \leq Z_j \leq x/A_j\}$$
$$= \frac{c}{x^{\alpha+\delta}} \mathrm{E} \sum |\tilde{A}_j|^{\alpha+\delta} \mathrm{E}_Z(|Z_j|^{\alpha+\delta} I\{-h(x/A_j) \leq Z_j \leq x/A_j\})$$

and once more using (4.8), we can bound, for large $x$, this expression by

$$\frac{c}{x^{\alpha+\delta}} \mathrm{E} \sum |\tilde{A}_j|^{\alpha+\delta} \left(\frac{x}{|\tilde{A}_j|}\right)^{\alpha+\delta} P(|Z| > x/|\tilde{A}_j|)$$
$$\leq c P(|Z| > x) \mathrm{E} \sum |\tilde{A}_j|^{\alpha-\varepsilon} \leq c P(|Z| > x) \mathrm{E} \left(\sum |\tilde{A}_j|^2\right)^{(\alpha-\varepsilon)/2}.$$

Therefore, by (3.9),

$$\lim_{n\to\infty} \limsup_{x\to\infty} \frac{\mathrm{II}_a}{P(|Z| > x)} = 0.$$

Next, it follows from independence and (4.12) that for some $0 < C < \infty$,

$$\mathrm{E}(Z_{j_1}^2 \ldots Z_{j_k}^2) \leq C \tag{4.13}$$

for all $(j_1,\ldots,j_k) \in D_k$. Therefore,

$$\mathrm{II}_b \leq \frac{c}{x^{\alpha+\delta}} \mathrm{E}\Bigg[\sum_{(j_1,\ldots,j_k)\in D_k} \tilde{A}_{j_1}^2 \ldots \tilde{A}_{j_k}^2\Bigg]^{(\alpha+\delta)/(2k)} \leq \frac{c}{x^{\alpha+\delta}} \mathrm{E}\left(\sum \tilde{A}_j^2\right)^{(\alpha+\delta)/2}$$

and so for $\delta$ small enough, by (3.9),

$$\lim_{x\to\infty} \frac{\mathrm{II}_b}{P(|Z| > x)} = 0.$$

This takes care of the term II and the argument for III is, as we have seen a number of times before, entirely similar.

Finally, let us consider the case $\alpha = 2m$ for some integer $m > 1$. This case is very similar to the case $\alpha > 2$ and not equal to an even integer, however, (4.12) does not hold ($k = m+1$ for small $\delta$ here). This does not make a difference as far as the term $\mathrm{II}_a$ above is concerned. For the term $\mathrm{II}_b$, we proceed as follows. Write $\hat{D}_k$ for the subset of $D_k$ where exactly $k-1$ of the indices are equal. Note that for $\delta$ small enough, the bound (4.13) still holds for all $(j_1,\ldots,j_k) \in D_k \setminus \hat{D}_k$. Then

$$\mathrm{II}_b \leq \frac{c}{x^{\alpha+\delta}} \mathrm{E}\left[\sum_{(j_1,\ldots,j_k) \in \hat{D}_k} \tilde{A}_{j_1}^2 \ldots \tilde{A}_{j_k}^2 Z_{j_1}^2 \ldots Z_{j_k}^2 \right.$$
$$\left. \times I\{-h(x/\tilde{A}_{j_i}) \leq Z_{j_i} \leq x/\tilde{A}_{j_i} \text{ for } i=1,\ldots,k\}\right]^{(\alpha+\delta)/(2k)}$$
$$+ \frac{c}{x^{\alpha+\delta}} \mathrm{E}\left[\sum_{(j_1,\ldots,j_k) \in D_k \setminus \hat{D}_k} \tilde{A}_{j_1}^2 \ldots \tilde{A}_{j_k}^2 Z_{j_1}^2 \ldots Z_{j_k}^2 \right.$$
$$\left. \times I\{-h(x/\tilde{A}_{j_i}) \leq Z_{j_i} \leq x/\tilde{A}_{j_i} \text{ for } i=1,\ldots,k\}\right]^{(\alpha+\delta)/(2k)}$$
$$:= \mathrm{II}_{b1} + \mathrm{II}_{b2}.$$

The term $\mathrm{II}_{b2}$ is treated in the same way as the term $\mathrm{II}_b$ with $\alpha > 2$ not being an even integer. Furthermore,

$$\mathrm{II}_{b1} \leq \frac{c}{x^{\alpha+\delta}} \mathrm{E}\left[\sum_{j_1}\sum_{j_2 \neq j_1} \tilde{A}_{j_1}^2 |\tilde{A}_{j_2}|^\alpha Z_{j_1}^2 |Z_{j_2}|^\alpha I\left\{-h\left(\frac{x}{\tilde{A}_{j_2}}\right) \leq Z_{j_2} \leq \frac{x}{\tilde{A}_{j_2}}\right\}\right]^{(\alpha+\delta)/(2k)}$$

$$\leq \frac{c}{x^{\alpha+\delta}} \mathrm{E}\left[\sum_{j_1}\sum_{j_2 \neq j_1} \tilde{A}_{j_1}^2 |\tilde{A}_{j_2}|^\alpha \mathrm{E}_Z\left(Z_{j_1}^2 |Z_{j_2}|^\alpha I\left\{-h\left(\frac{x}{\tilde{A}_{j_2}}\right) \leq Z_{j_2} \leq \frac{x}{\tilde{A}_{j_2}}\right\}\right)\right]^{(\alpha+\delta)/(2k)}$$

$$= \frac{c}{x^{\alpha+\delta}} \mathrm{E}\left[\sum_{j_1}\sum_{j_2 \neq j_1} \tilde{A}_{j_1}^2 |\tilde{A}_{j_2}|^\alpha \mathrm{E}_Z\left(|Z_{j_2}|^\alpha I\left\{-h\left(\frac{x}{\tilde{A}_{j_2}}\right) \leq Z_{j_2} \leq \frac{x}{\tilde{A}_{j_2}}\right\}\right)\right]^{(\alpha+\delta)/(2k)}.$$

By Karamata's theorem, $l(x) = \mathrm{E}|Z|^\alpha I\{|Z| \leq x\}$ is slowly varying at infinity and, as such, is bounded from above for large $x$ by $cx^\varepsilon$. Therefore, for large $x$,

$$\mathrm{II}_{b1} \leq \frac{c}{x^{\alpha+\delta}} \mathrm{E}\left[\sum_{j_1}\sum_{j_2 \neq j_1} \tilde{A}_{j_1}^2 |\tilde{A}_{j_2}|^\alpha \left(\frac{x}{|\tilde{A}_{j_2}|}\right)^\varepsilon\right]^{(\alpha+\delta)/(2k)}$$

$$= \frac{c}{x^{\alpha+\delta-\varepsilon(\alpha+\delta)/(2k)}} \mathrm{E}\left[\sum_{j_1} \tilde{A}_{j_1}^2 \sum_{j_2} |\tilde{A}_{j_2}|^{\alpha-\varepsilon}\right]^{(\alpha+\delta)/2k}$$

$$\leq \frac{c}{x^{\alpha+\delta-\varepsilon(\alpha+\delta)/(\alpha+2)}} \mathrm{E}\left(\sum \tilde{A}_j^2\right)^{(\alpha+2-\varepsilon)(\alpha+\delta)/(2(\alpha+2))}$$

and we see that for $\delta$, $\varepsilon$ small enough and $\varepsilon < \delta(\alpha+2)/(\alpha+\delta)$, the power of $x$ is larger than $\alpha$, so we may use (3.9) to obtain

$$\lim_{x \to \infty} \frac{\mathrm{II}_{b1}}{P(|Z| > x)} = 0$$

for all $n$. This completes the treatment of the term II in the case where $\alpha$ is an even integer greater than 2 and the term III is treated in the same way.

This proves the limit (4.2) in all cases and, hence, we have established the one-dimensional statement (4.1).

### 4.5. The multidimensional case

We now prove the general statement (3.10). For finite $n \geq 1$, Lemma 4.3 implies

$$\frac{P(u^{-1} \sum_{j \leq n} A_j Z_j \in \cdot)}{P(|Z| > u)} \xrightarrow{v} \mathrm{E}\left[\sum_{j \leq n} \mu \circ A_j^{-1}(\cdot)\right].$$

As in the one-dimensional case considered above, (3.10) will follow once we check that

$$\lim_{n \to \infty} \limsup_{x \to \infty} \frac{P(|\sum_{j > n} A_j Z_j| > x)}{P(|Z| > x)} = 0.$$

Because of the finite-dimensionality, it is enough to prove the corresponding statement for each one of the $d \times p$ elements of the random matrices. We will check (in the obvious notation) that

$$\lim_{n \to \infty} \limsup_{x \to \infty} \frac{P(|\sum_{j > n} A_j^{(11)} Z_j^{(1)}| > x)}{P(|Z| > x)} = 0. \tag{4.14}$$

Note that $|A_j^{(11)}| \leq \|A_j\|$ and so the sequence $(A_j^{(11)})$ in (4.14) satisfies the one-dimensional version of the assumptions (3.7)–(3.9). Only one thing prevents us from immediately applying the one-dimension statement (4.1) and it is the fact that the tail of $|Z^{(1)}|$ may happen to be strictly lighter than that of $\|Z\|$. To overcome this problem, let $(\tilde{Z}_j)$ be an independent copy of the sequence $(Z_j)$, also independent of the sequence $(A_j^{(11)})$, and $(\varepsilon_j)$ a sequence of i.i.d. Rademacher random variables, independent of the rest of the random variables. Then

$$P\left(\left|\sum_{j>n} A_j^{(11)} Z_j^{(1)}\right| > x\right) \leq P\left(\left|\sum_{j>n} A_j^{(11)}(Z_j^{(1)} + \varepsilon_j|\tilde{Z}_j|)\right| > x/2\right)$$
$$+ P\left(\left|\sum_{j>n} A_j^{(11)} \varepsilon_j|\tilde{Z}_j|\right| > x/2\right),$$

which allows us to apply (4.1) to each term above and prove (4.14). This completes the proof of the theorem.

*Remark 4.1.* The argument leading to (4.1) and (4.2) can also be used to show various modifications of these two statements. For example, in (4.2), one can allow $n$ and $x$ to go to infinity at the same time to obtain

$$\lim_{x \to \infty} \frac{P(\sum_{j>n(x)} A_j Z_j > x)}{P(|Z| > x)} = 0 \tag{4.15}$$

for any function $n(x) \to \infty$ as $x \to \infty$. Further, only values of the $Z_j$'s comparable to the level $x$ are relevant in (4.1), in the sense that

$$\lim_{\tau \to 0} \limsup_{x \to \infty} \frac{P(|\sum_{j=1}^{\infty} A_j Z_j 1(|Z_j| \leq \tau x)| > x)}{P(|Z| > x)} = 0. \tag{4.16}$$

**Lemma 4.1.** *Let $Z \in \mathrm{RV}(\alpha, \mu)$ be a random vector with $0 < \alpha < 1$. Then for each $y > 0$ $\varepsilon > 0$, and $c > 0$, there is some $x_0$ such that*

$$\frac{|\mathrm{E}[ZI\{|Z| \leq x/y\}]|}{xP(|Z| > x)} \leq c \max[y^{\alpha-1+\varepsilon}, y^{\alpha-1-\varepsilon}], \qquad x \geq x_0.$$

**Proof.** Karamata's theorem implies

$$\mathrm{E}[|Z|I\{|Z| \leq x\}] \sim \alpha(1-\alpha)^{-1} x P(|Z| > x)$$

and the claim follows from Potter's bound. $\square$

**Lemma 4.2.** *Let $Z \in \mathrm{RV}(\alpha, \mu)$ with $\alpha > 1$ a one-dimensional continuous random variable and such that $\mu$ assigns positive weights to both $(-\infty, 0)$ and $(0, \infty)$. Then there are numbers $K, C > 0$ and a function $h : [K, \infty) \to (0, \infty)$ satisfying $C^{-1} \leq h(x)/x \leq C$ and*

$$\int_x^{\infty} y P(Z \in \mathrm{d}y) = \int_{-\infty}^{-h(x)} |y| P(Z \in \mathrm{d}y)$$

*for all $x \geq K$.*

**Proof.** By the assumptions, $B = \int_{-\infty}^{0} |y| P(Z \in \mathrm{d}y) \in (0, \infty)$. The function $G(x) = \int_x^{\infty} y P(Z \in \mathrm{d}y)$ is continuous and decreases to zero. Let $K \geq 1$ be such that $G(x) \leq B/2$ for $x \geq K$ and define

$$h(x) = \inf\left\{t > 0 : \int_{-\infty}^{-t} |y| P(Z \in \mathrm{d}y) = G(x)\right\}.$$

It remains to prove the existence of a number $C$ in the statement. Let $g(x) = h(x)/x$. First, we show that $g(x)$ is bounded. By Karamata's theorem,

$$H(x) = \int_{-\infty}^{-h(x)} |y| P(Z \in \mathrm{d}y) \sim (1-w)h(x)^{1-\alpha} l(h(x))$$
$$= (1-w)(g(x)x)^{1-\alpha} l(g(x)x),$$
$$G(x) \sim w x^{1-\alpha} l(x),$$

where $l(x)$ is a slowly varying function. Supposing that $g(x) \to \infty$, then since $H(x) = G(x)$, we can find a constant $C$ and $\varepsilon \in (0, \alpha - 1)$ such that

$$1 = \frac{H(x)}{G(x)} \sim \frac{1-w}{w} g(x)^{1-\alpha} \frac{l(g(x)x)}{l(x)} \leq C g(x)^{1-\alpha} g(x)^{\varepsilon}.$$

But the right-hand side converges to 0, which is a contradiction. Hence, $g(x)$ cannot be unbounded. A similar argument shows that $g(x)$ must be bounded away from 0. □

The following is a slight generalization of Basrak *et al.* (2002), Proposition A.1.

**Lemma 4.3.** *Let $Z_1, \ldots, Z_n$ be i.i.d. random vectors in $\mathrm{RV}(\alpha, \mu)$ and $A_1, \ldots, A_n$ a sequence of random matrices such that for every $j = 1, \ldots, n$, $Z_j$ is independent of $\sigma(A_1, \ldots, A_j)$. Assume, further, that for some $\varepsilon > 0$, $\mathrm{E}\|A_j\|^{\alpha+\varepsilon} < \infty$ for $j = 1, \ldots, n$. Then*

$$\frac{P(u^{-1} \sum_{j=1}^n A_j Z_j \in \cdot)}{P(|Z| > u)} \xrightarrow{v} \mathrm{E}\left[\sum_{j=1}^n \mu \circ A_j^{-1}(\cdot)\right]. \tag{4.17}$$

**Proof.** Let $\hat{\mu}$ denote the measure in the right-hand side of (4.17) (as usual, the origin is not a part of the space where $\hat{\mu}$ lives, so any mass at the origin is simply lost). Let $B \in \mathbb{R}^d$ be a Borel set, bounded away from the origin, and assume that it is a $\hat{\mu}$-continuity set. We will use the notation $A^\varepsilon = \{x \in \mathbb{R}^d : d(x, A) < \varepsilon\}$ and $A_\varepsilon = \{x \in A : d(x, A^c) > \varepsilon\}$ for a set $A \in \mathbb{R}^d$ and $\varepsilon > 0$. Choose $\varepsilon > 0$ so small that $B^\varepsilon$ is still bounded away from the origin. We have

$$P\left(u^{-1} \sum_{j=1}^n A_j Z_j \in B\right) \leq P\left(\bigcup_{j=1}^n \{u^{-1} A_j Z_j \in \overline{B}^\varepsilon\}\right)$$
$$+ P\left(\bigcup_{j_1=1}^n \bigcup_{j_2=j_1+1}^n \{u^{-1} |A_{j_1} Z_{j_1}| > \varepsilon/n, u^{-1} |A_{j_2} Z_{j_2}| > \varepsilon/n\}\right).$$

By the Portmanteau theorem,

$$\limsup_{u \to \infty} \frac{P(\bigcup_{j=1}^n \{u^{-1} A_j Z_j \in \overline{B}_\varepsilon\})}{P(|Z| > u)} \leq \hat{\mu}(\overline{B}^\varepsilon).$$

On the other hand, by Basrak *et al.* (2002), Proposition A.1 for every $j_1 < j_2$ and $M > 0$,

$$\limsup_{u \to \infty} \frac{P(u^{-1}|A_{j_1}Z_{j_1}| > \varepsilon, u^{-1}|A_{j_2}Z_{j_2}| > \varepsilon)}{P(|Z| > u)}$$
$$\leq \limsup_{u \to \infty} \frac{P(|A_{j_1}Z_{j_1}| > M, |A_{j_2}Z_{j_2}| > u\varepsilon)}{P(|Z| > u)}$$
$$= \varepsilon^{-\alpha} \mathrm{E}(|A_{j_2}|^\alpha I\{|A_{j_1}Z_{j_1}| > M\})$$

and letting $M \to \infty$, we obtain

$$\lim_{u \to \infty} \frac{P(u^{-1}|A_{j_1}Z_{j_1}| > \varepsilon, u^{-1}|A_{j_2}Z_{j_2}| > \varepsilon)}{P(|Z| > u)} = 0.$$

Using regular variation, we conclude that

$$\limsup_{u \to \infty} \frac{P(u^{-1} \sum_{j=1}^n A_j Z_j \in B)}{P(|Z| > u)} \leq \hat{\mu}(\overline{B^\varepsilon}).$$

Letting $\varepsilon \to 0$ and using the fact that $B$ is a $\hat{\mu}$-continuity set, we obtain

$$\limsup_{u \to \infty} \frac{P(u^{-1} \sum_{j=1}^n A_j Z_j \in B)}{P(|Z| > u)} \leq \hat{\mu}(B).$$

A matching lower bound follows in a similar way using the relation

$$P\left(u^{-1} \sum_{j=1}^n A_j Z_j \in B\right) \geq P\left(\bigcup_{j=1}^n \{u^{-1} A_j Z_j \in B_\varepsilon\}\right)$$
$$- P\left(\bigcup_{j_1=1}^n \bigcup_{j_2=j_1+1}^n \{u^{-1}|A_{j_1}Z_{j_1}| > \theta, u^{-1}|A_{j_2}Z_{j_2}| > \theta\}\right),$$

where $\theta = \min\{\varepsilon/n, \inf_{x \in B} \|x\|\}$. □

## Acknowledgements

The authors gratefully acknowledge the comments made by the anonymous referees. Henrik Hult's research was supported in part by the Swedish Research Council and the Sweden America Foundation. Gennady Samorodnitsky's research is supported in part by NSF Grant DMS-03-03493 and NSA Grant H98230-06-1-0069 at Cornell University.